\documentstyle[amstex,12pt]{article}

\begin{document}
\catcode`\@=\active
\catcode`\@=11

\def\@eqnnum{\hbox to .01pt{}\rlap{\bf \hskip -\displaywidth(\theequation)}}

\newcommand{\C}{\mbox{$\Bbb C$}}
\newcommand{\Pj}{\mbox{$\Bbb P$}}
\newcommand{\Q}{\mbox{$\Bbb Q$}}
\newcommand{\Z}{\mbox{$\Bbb Z$}}
\newcommand{\Gammab}{\mbox{$\bar \Gamma$}}
\newcommand{\Deltab}{\mbox{$\bar \Delta$}}
\newcommand{\fp}{\mbox{     $\Box$}}
\newcommand{\lrou}{\mbox{$\ulcorner\!$}}
\newcommand{\rrou}{\mbox{$\urcorner$}}
\newcommand{\Ext}{\mbox{Ext}}

\begin{center}
 {\bf Bogomolov Instability and Kawamata-Viehweg Vanishing} \\
 \vskip.4cm
 {\sc Guillermo Fern\'andez del Busto}\footnote{Partially supported by UNAM} \\
\end{center}

\vskip.8cm
\noindent
{\bf Introduction}

\vskip.4cm
The purpose of this note is to show how the Kawamata-Viehweg vanishing theorem
for fractional divisors leads to a quick new proof of Bogomolov's instability
theorem for rank two vector bundles on an algebraic surface.

\vskip.4cm
Let $X$ be a smooth complex projective surface and let $\cal E$ be a rank two
holomorphic vector bundle on $X$. Bogomolov's theorem states that if
$c_1(\cal E)^2 >4c_2(\cal E)$, then $\cal E$ satisfies a strong instability
condition, which roughly speaking means that $\cal E$ contains an
exceptionally positive rank one subbundle.
Bogomolov's original proof \cite{B,Rd} revolved around a beautiful argument
with geometric invariant theory. Another proof---using characteristic $p$
techniques---was given by Gieseker \cite{G}, and Miyaoka \cite{M} subsequently
found a simple way to reduce the question to some restriction theorems of
Mumford, Mehta and Ramanathan \cite{MR}. Each of these proofs also yields
an analogous assertion for higher rank bundles.

\vskip.4cm
It is well understood that the rank two case of Bogomolov's result can be
used to prove various sorts of vanishing theorems on the surface $X$. For
example, Mumford \cite{Rd} used this approach to give a quick proof of
Ramanujam's form of the Kodaira vanishing theorem. Somewhat later, Reider
\cite{Rdr} realized that similar techniques yield criteria for the vanishing of
groups of the form $H^1(X, \cal O_X(K_X+L) \otimes \cal I_Z)$, where
$\cal I_Z$ is the ideal sheaf of a set of points $Z \subseteq X$.
He deduces thereby his celebrated theorem on freeness and very ampleness of
adjoint linear series on $X$. A cohomological approach to these questions,
based on Miyaoka's vanishing theorem for Zariski decompositions, was given
by Sakai \cite{S}.

\vskip.4cm
More recently, in connection with three dimensional analogues of Reider's
results, Ein and Lazarsfeld \cite{EL} found that one can prove (special cases
of) Reider's results using the Kawamata-Viehweg vanishing theorem. The argument
draws on the cohomological techniques pioneered by Kawamata, Koll\'ar, Reid and
Shokurov in connection with the minimal model program (c.f. \cite{CKM} or
\cite{KMM}). The question then arises whether similar techniques can be used to
deduce the full theorem of Bogomolov. We  complete this circle of ideas by
showing that this is indeed the case, and obtain a transparent new proof
of Bogomolov's theorem.

\vskip.4cm
The idea of the proof is very simple. Given the rank two bundle $\cal E$, after
twisting by a sufficiently ample line bundle we may assume that $\cal E$ has
lots of sections.
The Koszul complex associated to a general section $s$ then expresses $\cal E$
as
sitting in an exact sequence
$$0\rightarrow \cal O_X \stackrel{s}\rightarrow \cal E \rightarrow \cal O_X(L)
    \otimes \cal I_Z \rightarrow 0 \phantom{.},$$
where $L$ is ample, and $Z$ consists of distinct points. As in the proof of
Reider's theorem in~\cite{EL}, the positivity of
$c_1^2 (\cal E)-4c_2(\cal E)$ implies the existence of a divisor $D$ in $|nL|$
with high
multiplicity on $Z$. If the singular points of $D$ are (close to being)
isolated, then one would obtain a vanishing which contradicts the local
freeness of $\cal E$. Consequently, $D$ must have some special components
appearing with high multiplicity. We use these distinguished components to
construct a divisor $\Gamma \subseteq D$ through $Z$. The Kawamata-Viehweg
vanishing theorem implies then that the inclusion
$\cal O_X(L-\Gamma) \hookrightarrow \cal O_X(L) \otimes \cal I_Z$ lifts to an
embedding of $\cal O_X(L-\Gamma)$ into $\cal E$, and we argue that
(the saturation of) $\cal O_X(L-\Gamma)\subseteq \cal E$ is a destabilizing
subsheaf.
A new feature in this approach is that we need to be somewhat careful in our
choice of $Z$ and $D$. We use a monodromy argument to show that for a general
choice
of $s$ and $D$, $D$ will have the same multiplicity at each point of $Z=Z(s)$.

\vskip.4cm
This note is part of my Ph.D. thesis at UCLA, and I would like to use this
opportunity to express my gratitude to Rob Lazarsfeld, not just for all the
knowledge he has shared
with me, but for his patience, guidance and encouragement. I would further like
to thank Lawrence Ein for valuable discussions and suggestions.

\vskip.8cm
\noindent
{\bf \S 1. The Uniform Multiplicity Property and Beginning of Proof}.

\vskip.4cm
Consider a rank 2 vector bundle $\cal E$ on a projective and nonsingular
surface $X$;
recall that $\cal E$ is said to be {\em (Bogomolov)-unstable} if there exists
an exact
sequence
$$0\rightarrow \cal O_X(M) \rightarrow \cal E
     \rightarrow \cal O_X(N) \otimes \cal I_{\chi} \rightarrow 0 \phantom{.},$$
with $M$ and $N$ divisors, and $\chi \subseteq X$ a 0-dimensional subscheme
with sheaf of ideals $\cal I_{\chi}$, such that
$M-N$ is in the positive cone of the N\'eron-Severi group, i.e.
$(M-N)^2>0$ and $(M-N)\cdot H>0$ for any ample divisor $H$ on $X$.
Bogomolov's theorem then asserts that an equivalent condition for $\cal E$ to
be
unstable is that $c_1(\cal E)^2>4c_2(\cal E)$.

\vskip.4cm
Since Bogomolov's theorem is invariant under tensor product by line bundles, we
may
assume that $\cal E$ is globally generated. Let $s\in H^0(X,\cal E)$ be a
general section; the Kozsul resolution associated to $s$ defines $\cal E$ as an
extension
\begin{equation}
   0\rightarrow \cal O_X \stackrel{s}\rightarrow \cal E \rightarrow \cal O_X(L)
     \otimes \cal I_Z \rightarrow 0 \label{a}
\end{equation}
with $Z=Z(s)$ the 0-scheme of $s$; we may suppose that $L=c_1(\cal E)$ is ample
and
that $Z$ consists of $c_2(\cal E)$ distinct points.

\vskip.4cm
The local freeness of $\cal E$ imposes some conditions on the finite set $Z$
(c.f.~\cite{GH} or~\cite{OSS}).
In fact, there is an element $e\in \Ext^1(\cal O_X(L)\otimes \cal I_Z, \cal
O_X)$
corresponding to the extension class~(\ref{a}). If $Z'\subseteq Z$ is any
proper
subset (possibly empty) the local freeness of $\cal E$ then implies that $e$ is
not in
the image of
$\Ext^1(\cal O_X(L)\otimes \cal I_{Z'},\cal O_X)\rightarrow
  \Ext^1(\cal O_X(L)\otimes \cal I_Z,\cal O_X)$.
Since $\Ext^1(\cal O_X(L)\otimes\cal I_Z,\cal O_X)$ is Grothendieck dual to
$H^1(X,\cal O_X(K_X+L)\otimes\cal I_Z)$, this is equivalent of the
non-surjectivity of
the evaluation map
\begin{equation}
	  \label{b}
   H^0(X,\cal O_X(K_X+L)\otimes \cal I_{Z'}) \longrightarrow
   H^0(Z-Z',\cal O_{Z-Z'}(K_X+L)) \phantom{.},
\end{equation}
i.e. $Z-Z'$ cannot impose independent conditions on
$|\cal O_X(K_X+L)\otimes \cal I_{Z'}|$.

\vskip.4cm
Suppose now that $c_1(\cal E)^2 > 4c_2(\cal E)$. As in the proof of Reider's
theorem in \cite{EL}, for $n\!\gg\! 0$ there is a divisor $D$ in $|nL|$ with
$mult_z(D)\geq 2n+1$
for any $z\in Z$. (In brief, by Riemann-Roch
$h^0(X,\cal O_X(nL))\sim \frac {n^2}2 L^2$, whereas it is
$\binom {2n+2}{2} \deg(Z) \sim \frac {n^2}2 4\deg(Z)$ conditions on a divisor
to
have multiplicity at least $2n+1$ at each point of $Z$.)

\vskip.4cm
The proof will involve analysis of $D$ and $s$, but we will need to have some
control on the geometry of $D$. In fact, we need to show that
if the section $s\! \in \! H^0(X,\cal E)$ and the divisor $D\! \in \! |nL|$ are
sufficiently
general, then $D$ and $Z=Z(s)$ satisfy the following

\vskip.4cm
\noindent
{\sc Uniform Multiplicity Property (UMP).} For any rational number $\delta >0$,
the multiplicity of $[\delta D]$ is the same at every point of $Z$.

\vskip.4cm
For the proof,  let $\Pj (=\Pj H^0(\cal E))$ be the projective space
parametrizing
sections of $\cal E$. Let
$$X\times \Pj \supseteq \cal Z =\{ (x,s)\colon s(x)=0\}$$
be the {\em universal 0-scheme}, with projections $p$ and $q$.
Since $\cal E$ is globally generated,
$\cal Z\rightarrow X$ is a projective bundle, hence $\cal Z$ is
irreducible. Note that $\cal Z$ is generically finite over $\Pj$ and that
$\dim \cal Z =\dim \Pj$. Now let $\cal D\subseteq X\times \Pj $ be an
effective divisor flat over $\Pj$: we'll write $D_s$ for the fibre of $\cal D$
over $s$. The situation is as follows:
$$\begin{array}{rcc}
   &&\cal D \\
   &&\cap \\
  \cal Z\phantom{..}&\subset&X\times \Pj \\
  \searrow \! &&\swarrow \hfill \\
   &\Pj& \\
  \end{array}$$
The main point is the following

\vskip.4cm
\noindent
{\bf Lemma.} Let $k$ be a positive integer, and suppose that for a general
$s\in \Pj$,
there is a point $x\in Z(s)$ such that $mult_x(D_s)\geq k$. Then for a
general $s\in \Pj$, $mult_x(D_s)\geq k$ for {\em every} $x\in Z(s)$.

\vskip.4cm
\noindent
{\it Proof.} Let $\cal Z\supseteq \cal Z_k =
                   \{(x,s)\colon mult_x(D_s)\geq k\}$.
By hypothesis $\cal Z_k$ dominates $\Pj$, hence
$$\dim \cal Z_k=\dim \Pj =\dim \cal Z$$
and since $\cal Z$ is irreducible, we have that $\cal Z_k=\cal Z$. \fp

\vskip.4cm
For the {\sc (UMP)}, consider the sheaf
$\cal F=p^*\cal O_X(nL)\otimes \cal I_{\cal Z}^{2n+1}$. Since
$c_1(\cal E)^2>4c_2(\cal E)$, $q_*\cal F$ has positive rank. Hence
$H^0(X\times \Pj,q^*\cal O_{\Pj}(H)\otimes \cal F)\not= 0$
for sufficiently positive $H$. Let $\cal D$ be the corresponding divisor
of a section of $q^*\cal O_{\Pj}(H)\otimes \cal F$. Now apply the previous
lemma to the divisor $[\delta \cal D]$.

\vskip.4cm
The basic tool for our cohomological approach to Bogomolov instability is
the vanishing theorem of Kawamata-Viehweg. Sakai noticed (c.f.~\cite[1.1]{EL})
that on a surface $X$, the vanishing of $\Q$-divisors holds without any
normal crossing assumption:

\vskip.4cm
\noindent
{\bf Theorem(Kawamata-Viehweg)}. Let $X$ be a nonsingular projective surface,
and let $M$ be any big and nef $\Q$-divisor on $X$. Then
\begin{equation}
\label{j}
  H^i(X,\cal O_X(K_X+\lrou M\rrou ))=0
  \phantom{.}\mbox{for} \phantom{.} i>0 \phantom{.}.
\end{equation}

\vskip.8cm
\noindent
{\bf \S 2. End of Proof.}

\vskip.4cm
With notation as in \S 1 now we fix a general section $s\in \Pj $, and write
$Z=Z(s)$
and $D=D_s$. We assume that the {\sc (UMP)} holds for $Z$ and $D$. The
{\sc (UMP)} implies that $D$ has the same multiplicity at every point $z\in Z$;
denote by $m$ such multiplicity. Let
$$D=\sum d_jD_j+F\phantom{.},$$
where $D_j$ are the components of $D$ intersecting $Z$, and $F$ consists of the
components of $D$ disjoint from $Z$. Let $d=\max \{ d_j \}$ and consider the
divisor
$$D_0=[\sum \frac{d_j}d D_j]\phantom{.}.$$
Then $D_0\subseteq X$ is reduced and contains (by construction) at least one
point of
$Z$, hence $Z\subseteq D_0$ because of the {\sc (UMP)}.

\vskip.4cm
Notice next that $2d>m$. In fact, if $2d\leq m$ then as in the proof of
Reider's
theorem in \cite{EL}, let $f\colon Y\rightarrow X$ be the blow-up of $X$ at $Z$
and
let $E_i\subseteq Y$ be the exceptional divisors. The $\Q$-divisor
$f^*(L-\frac {2\delta_1}m D)- \delta_2 \sum E_i\equiv
 \delta_1f^*(L-\frac 2m D)+(1-\delta_1)f^*L-\delta_2 \sum E_i$,
with  $0\! <\! \delta_2 \! \ll \!  \delta_1 \!< \!1$ and
$\delta_2  \! \geq \!2(1-\delta_1)$, is nef and big, and applying the
Kawamata-Viehweg vanishing theorem~(\ref{j}), it follows that
$$H^1(X,\cal O_X(K_X+L-[\frac 2m F])\otimes \cal I_Z)=0\phantom{.},$$
which is absurd by~(\ref{b}).

\vskip.4cm
Assume henceforth that $2d>m$. In particular, $D_0$ is nonsingular at each
point of
$Z$ and each point of $Z$ is in exactly one component of $D_0$.
The $\Q$-divisor $L-\frac 1d D$ is ample, so $H^1(X,\cal O_X(K_X+L-[\frac 1d
D]))=0$. Therefore there exists a minimal subdivisor $\Delta \subseteq [\frac
1d F]$ such that
\begin{equation}
\label{c}
   H^1(X,\cal O_X(K_X+L-D_0-\Delta))=0\phantom{.},
\end{equation}
in the sense that if~(\ref{c}) holds for $\Delta$, and if $\Delta '$ is a
component of $\Delta$ (whenever $\Delta \not= 0$) then
$H^1(X,\cal O_X(K_X+L-D_0-\Delta +\Delta ')) \not= 0$. Let
$$\Gamma =D_0+\Delta \phantom{.}.$$
Then~(\ref{c}) implies that $\Ext^1(\cal O_X(L-\Gamma),\cal O_X)=0$, and in
consequence there is an injection $\cal O_X(L-\Gamma) \hookrightarrow \cal E$
with
$$\begin{array}{ccccccccc}
  &&&&&&\cal O_X(L-\Gamma)&& \\
  &&&&&\swarrow&\downarrow&& \\
  0&\longrightarrow&\cal O_X&
      \longrightarrow&\cal E &
      \longrightarrow&\cal O_X(L)\otimes\cal I_Z&
      \longrightarrow &0 \phantom{.}. \\
   \end{array}$$
The purpose now is to show that (the saturation of) $\cal O_X(L-\Gamma)$ is the
destabilizing subsheaf of $\cal E$. For this end, we need to establish some
inequalities.

\vskip.4cm
Note first that if $\Delta'$ is a component of $\Delta$, then
\begin{equation}
   \label{d}
   (L-\Gamma)\cdot \Delta' \leq 0\phantom{.}.
\end{equation}
In fact, the exact sequence
$$0\rightarrow \cal O_X(K_X+L-\Gamma) \rightarrow
     \cal O_X(K_X+L-\Gamma+\Delta')\rightarrow
     \cal O_{\Delta'}(K_{\Delta'}+L-\Gamma)\rightarrow 0\phantom{.},$$
and $H^1(X,\cal O_X(K_X+L-\Gamma+\Delta'))\not=0$ imply that
$H^1(\Delta',\cal O_{\Delta'}(K_{\Delta'}+(L-\Gamma)))\not=0$,
whence~(\ref{d}).
We claim next that
\begin{equation}
  \label{e}
  (L-\Gamma)\cdot D_0 \leq \deg(Z) \phantom{.}.
\end{equation}
In fact, since each point of $Z$ is in exactly one component of $D_0$ and each
component of $D_0$ contains at least one point of $Z$, then~(\ref{e}) will
follow if
we show that if $D'$ is a component of $D_0$ containing a subset $Z'\subseteq
Z$, then $(L-\Gamma)\cdot D' \leq \deg(Z')$.
Suppose to the contrary that $(L-\Gamma)\cdot D' >\deg(Z')$. Then any
$\deg(Z')$
nonsingular points on $D'$ impose independent conditions on the linear series
$|\cal O_{D'}(K_{D'}+L-\Gamma)|$.
On the other hand, it follows from~(\ref{c}) that
$$H^0(X,\cal O_X(K_X+L-\Gamma+D'))\longrightarrow
    H^0(D', \cal O_{D'}(K_{D'}+L-\Gamma))$$
is surjective. Since $Z-Z'\subseteq \mbox{Supp}(\Gamma -D')$ we have then that
$Z'$ imposes independent conditions on $H^0(X,\cal O_X(K_X+L)\otimes \cal
I_{Z-Z'})$.
But this contradicts~(\ref{b}), so~(\ref{e})
is established. Combining~(\ref{d}) and~(\ref{e}) yields
\begin{equation}
   \label{f}
   (L-\Gamma)\cdot \Gamma \leq \deg(Z) \phantom{.}.
\end{equation}

\vskip.4cm
We now assert that
\begin{equation}
   \label{g}
   \deg(Z)>\frac 12 L\cdot \Gamma \phantom{.}.
\end{equation}
For this, note that
$$L-\Gamma \equiv (1-\frac nd )L+(\sum \frac{d_j}d D_j-D_0)
                                                      +(\frac 1d
F-\Delta)\phantom{.}.$$
We assert first that
\begin{equation}
   \label{h}
   (\frac 1d F-\Delta)\cdot\Gamma\geq 0 \phantom{.}.
\end{equation}
In fact, let $\frac 1d F-\Delta=\Delta_1+\Delta_2$, with $\Delta_1$ an
effective $\Q$-divisor all of whose components are components of $\Delta$, and
with $\Delta_2$ and $\Delta$ having no common components. Then
$\Delta_2\cdot \Gamma\geq 0$, hence
$(\frac1dF-\Delta)\cdot\Gamma\geq\Delta_1\cdot\Gamma$ and
$\Delta_1\cdot\Gamma\geq\Delta_1\cdot L\geq 0$ because of~(\ref{d}). Recalling
that $\sum \frac{d_j}d D_j -D_0$ and $\Delta$ have no common components,
observe next:
$$\begin{array}{rcl}
        (\sum\frac{d_j}d D_j-D_0)\cdot\Gamma&\geq&
             (\sum\frac{d_j}d D_j-D_0)\cdot D_0 \\
        &\geq&\frac 1d \displaystyle{\sum_{z\in Z}}(mult_z(D)-d) \\
        &\geq&\deg(Z) (\frac md-1) \phantom{.}. \\
\end{array}$$
Combining~(\ref{f}) and~(\ref{h}) we then conclude that
$$\deg(Z)\geq (L-\Gamma)\cdot \Gamma
               \geq (1-\frac nd )L\cdot \Gamma + \deg(Z) (\frac md-1)$$
and since $2d\! >\! m\! >\! 2n$,~(\ref{g}) follows.

\vskip.4cm
Now consider the subsheaf $\cal O_X(L-\Gamma) \hookrightarrow \cal E$. The
saturation of this subsheaf is of the form $\cal O_X(L-\Gammab)$, for a
subdivisor
$\Gammab \subseteq \Gamma$. We claim that $\cal O_X(L-\Gammab)$ is a
destabilizing subsheaf of $\cal E$; this is equivalent to the inequalities
$$\begin{array}{rcl}
    (L-2\Gammab)^2&>&0 \\
    (L-2\Gammab)\cdot L&>&0 \phantom{.}. \\
   \end{array}$$
Now the first of this inequalities is a consequence of the hypothesis
$c_1^2(\cal E)>4c_2(\cal E)$ by computing $c_2(\cal E)$ from the exact sequence
$0\rightarrow \cal O_X(L-\Gammab) \rightarrow \cal E \rightarrow
     \cal O_X(\Gammab) \otimes \cal I_{Z '}\rightarrow 0 \phantom{.}.$

\vskip.4cm
\noindent
As for the second inequality, it is a consequence of~(\ref{g}):
$$L^2>4\deg(Z)>2L\cdot \Gamma \geq 2L\cdot \Gammab \phantom{.}.$$

\vskip.4cm
This completes the proof of Bogomolov's theorem.

\vskip.6cm
\noindent
{\it Remark.} A similar argument as in the proof of Bogomolov's theorem
completes
the proof of Reider's theorem in~\cite[1.4]{EL}.

\vskip.8cm
\begin{flushleft}
Department of Mathematics \\
University of California, Los Angeles \\
Los Angeles, CA 90024-1555 \\
e-mail: gonzalez@math.ucla.edu \\
\end{flushleft}


\begin{thebibliography}{KMM}
\bibitem[{\bf B}]{B} F. A. Bogomolov, {\em Holomorphic tensors and vector
bundles on
   projective varieties}, Math USSR Izvestija {\bf 13} (1979), {\em no. 3},
pp.499-555.
\bibitem[{\bf CKM}]{CKM} H. Clemens, J. Koll\'ar and S. Mori, {\em Higher
dimensional
   Complex Geometry}, Asterisque {\bf 166}, 1988.
\bibitem[{\bf EL}]{EL} L. Ein and R. Lazarsfeld, {\em Global generation of
   pluricanonical and adjoint linear series on smooth projective threefolds},
   preprint.
\bibitem[{\bf G}]{G} D. Gieseker, {\em On a theorem of Bogomolov on Chern
classes
   of stable bundles}, Amer. J. Math. {\bf 101} (1979), {\em no. 1}, pp. 77-85.
\bibitem[{\bf GH}]{GH} P. Griffiths and J. Harris, {\em Residues and
zero-cycles on
   algebraic varieties}, Ann. Math. {\bf 108} (1978), pp. 461-505.
\bibitem[{\bf KMM}]{KMM} Y. Kawamata, K. Matsuda and K. Matsuki, {\em
Introduction to
   the minimal model problem}, in {\em Algebraic Geometry, Sendai}, Adv.
Studies
   in Pure Math {\bf 10}, Kinokuniya-North-Holland (1987), pp. 283-360.
\bibitem[{\bf M}]{M} Y. Miyaoka, {\em The Chern classes and Kodaira dimension
of
   a minimal variety}, in {\em Algebraic Geometry, Sendai}, Adv. Studies
   in Pure Math {\bf 10}, Kinokuniya-North-Holland (1987), pp. 449-476.
\bibitem[{\bf MR}]{MR} V. Mehta and A. Ramanathan, {\em Semi-stable sheaves on
   projective varieties and their restriction to curves}, Math. Ann. {\bf 258}
   (1982), pp. 213-224.
\bibitem[{\bf OSS}]{OSS} C. Okoneck, H. Spindler and M. Schneider, {\em Vector
bundles
   on Complex Projective Spaces}, Progress in Math. {\bf 3} (1980),
Birkh\"{a}user.
\bibitem[{\bf Rd}]{Rd} M. Reid, {\em Bogomolov's theorem $c_1^2\leq 4c_2$}, in
   {\em Proc. Intern. Symp. on Alg. Geom.}, Kyoto (1977), pp. 633-642.
\bibitem[{\bf Rdr}]{Rdr} I. Reider, {\em Vector bundles of rank 2 and linear
systems
   on algebraic surfaces}, Ann. Math. {\bf 127} (1988), pp. 309-316.
\bibitem[{\bf S}]{S} F. Sakai, {\em Reider-Serrano's method on normal
surfaces}, in
   {\em Algebraic Geometry: Proceedings, L'Aquila 1988}, Lect. Notes in Math.
   {\bf 1417} (1990), pp. 301-319.
\end{thebibliography}
\end{document}